\documentclass[centertags,leqno]{article}

\usepackage{amsmath}
\usepackage{amssymb}
\usepackage{latexsym}
\usepackage{amsxtra}
\usepackage{amscd}
\usepackage{theorem}

\usepackage{curves}
\usepackage{epic}
\usepackage{eepic}

\theoremstyle{plain}

{\theorembodyfont{\slshape}

        \newtheorem{thm}{Theorem}[section]
        \newtheorem{cor}[thm]{Corollary}
        \newtheorem{lem}[thm]{Lemma}
        \newtheorem{prop}[thm]{Proposition}
        
}

{\theorembodyfont{\rmfamily}

        \newtheorem{rem}[thm]{Remark}
        
        \newtheorem{ass}[thm]{Assumption}
        
        \newtheorem{convention}[thm]{Convention}

}

\renewcommand{\em}{\sl}

\newcommand{\proof}{{\bf Proof:\ }}
\newcommand{\Endproof}{\hspace*{\fill} $\Box$ \vspace{1ex} \noindent }

\makeatletter
\renewcommand{\subsection}{\@startsection{subsection}{2}%
        {\z@}{-3.25ex plus -1ex minus-.2ex}{-1em}{\bf}}
\makeatother
\renewcommand{\em}{\sl}

\makeatletter
\renewcommand{\subsection}{\@startsection{subsection}{2}%
{\z@}{-3.25ex plus -1ex minus-.2ex}{-1em}{\bf}} \makeatother

\newcommand{\PP}{\mathbb{P}}
\newcommand{\ZZ}{\mathbb{Z}}
\newcommand{\CC}{\mathbb{C}}
\newcommand{\RR}{\mathbb{R}}
\newcommand{\QQ}{\mathbb{Q}}

\renewcommand{\AA}{\mathbb{A}}

\newcommand{\E}{\mathcal{E}}
\newcommand{\V}{\mathcal{V}}
\newcommand{\W}{\mathcal{W}}

\newcommand{\OO}{\mathcal{O}}

\newcommand{\VB}{\bar{\V}}

\newcommand{\g}{{\bf g}}

\newcommand{\GL}{{\rm GL}}

\newcommand{\PGL}{{\rm PGL}}

\newcommand{\supp}{{\rm supp}}

\newcommand{\pib}{\bar{\pi}}

\newcommand{\para}{_p}

\newcommand{\inj}{\hookrightarrow}
\newcommand{\To}{\;\longrightarrow\;}
\newcommand{\iso}{\stackrel{\sim}{\to}}
\newcommand{\liso}{\;\stackrel{\sim}{\longrightarrow}\;}

\newcommand{\gen}[1]{\mathopen\langle#1\mathclose\rangle}

\newcommand{\dR}{_{\rm dR}}

\title{Variation of parabolic cohomology and Poincar\'e duality}

\author{
   Michael Dettweiler \\ Universit\"at Heidelberg
   \and
   Stefan Wewers
 \\ Universit\"at Bonn}

\date{}

\begin{document}

\maketitle

\begin{abstract}
 We continue our study of the variation of parabolic cohomology
(\cite{dw04}) and derive an exact formula for the underlying Poincar\'e 
duality. As an illustration of our methods, we compute the monodromy of the 
Picard-Euler system and its invariant Hermitian form, reproving a
classical theorem of Picard.
\end{abstract}


\section*{Introduction}

Let $x_1,\ldots,x_r$ be pairwise distinct points on the Riemann sphere
$\PP^1(\CC)$ and set $U:=\PP^1(\CC)-\{x_1,\ldots,x_r\}$. The
Riemann--Hilbert correspondence \cite{DeligneED} is an equivalence
between the category of ordinary differential equations with
polynomial coefficients and at most regular singularieties at the
points $x_i$ and the category of local systems of $\CC$-vectorspaces
on $U$. The latter are essentially given by an $r$-tuple of matrices
$g_1,\ldots,g_r\in\GL_n(\CC)$ satisfying the relation
$\prod_ig_i=1$. The Riemann--Hilbert correspondence associates to a
differential equation the tuple $(g_i)$, where $g_i$ is the monodromy
of a full set of solutions at the singular point $x_i$.

In \cite{dw04} the authors investigated the following
situation. Suppose that the set of points
$\{x_1,\ldots,x_r\}\subset\PP^1(\CC)$ and a local system $\V$ with
singularities at the $x_i$ depend
on a parameter $s$ which varies over the points of a complex manifold
$S$. More precisely, we consider a relative divisor $D\subset\PP^1_S$
of degree $r$ such that for all $s\in S$ the fibre
$D_s\subset\PP^1(\CC)$ consists of $r$ distinct points. Let
$U:=\PP^1_S-D$ denote the complement and let $\V$ be a local system on
$U$. We call $\V$ a {\em variation of local systems} over the base
space $S$. The {\em parabolic cohomology} of the variation $\V$ is the
local system on $S$
\[
   \W \;:=\; R^1\pi_*(j_*\V),
\]
where $j:U\inj\PP^1_S$ denotes the natural injection and
$\pi:\PP^1_S\to S$ the natural projection. The fibre of $\W$ at a
point $s_0\in S$ is the parabolic cohomology of the local system
$\V_0$, the restriction of $\V$ to the fibre
$U_0=U\cap\pi^{-1}(s_0)$. 

A special case of this construction is the {\em middle convolution
functor} defined by Katz \cite{KatzRLS}. Here $S=U_0$ and so this
functor transforms one local system $\V_0$ on $S$ into another one, $\W$.
Katz shows that all rigid local systems on $S$ arise from
one-dimensional systems by successive application of middle
convolution. This was further investigated by Dettweiler and Reiter
\cite{DettwReiterMC}. Another special case are the generalized
hypergeometric systems studied by Lauricella \cite{Lauricella}, Terada
\cite{Terada} and Deligne--Mostow \cite{DeligneMostow}. Here $S$ is
the set of ordered tuples of pairwise distinct points on $\PP^1(\CC)$
of the form $s=(0,1,\infty,x_4,\ldots,x_r)$ and $\V$ is a
one-dimensional system on $\PP^1_S$ with regular singularities at the (moving)
points $0,1,\infty,x_4,\ldots,x_r$. In \cite{dw04} we gave another
example where $S$ is a $17$-punctured Riemann sphere and the local system
$\V$ has finite monodromy. The resulting local system $\W$ on $S$ does
not have finite monodromy and is highly non-rigid. Still, by the
comparison theorem between singular and \'etale cohomology, $\W$ gives
rise to $\ell$-adic Galois representations, with interesting
applications to the regular inverse Galois problem. 

In all these examples, it is a significant fact that the monodromy of
the local system $\W$ (i.e.\ the action of $\pi_1(S)$ on a fibre of
$\W$) can be computed explicitly, i.e.\ one can write down matrices
$g_1,\ldots,g_r\in\GL_n$ which are the images of certain generators
$\alpha_1,\ldots,\alpha_r$ of $\pi_1(S)$. In the case of the middle
convolution this was discovered by Dettweiler--Reiter
\cite{DettwReiterKatz} and V\"olklein \cite{Voelklein01}. In
\cite{dw04} it is extended to the more general situation sketched
above. In all earlier papers, the computation of the monodromy is
either not explicit (like in \cite{KatzRLS}) or uses ad hoc
methods. In contrast, the method presented in \cite{dw04} is very
general and can easily be implemented on a computer.

It is one matter to compute the monodromy of $\W$ explicitly (i.e.\ to
compute the matrices $g_i$) and another matter to determine its image
(i.e.\ the group generated by the $g_i$). In many cases the image of
monodromy is contained in a proper algebraic subgroup of $\GL_n$,
because $\W$ carries an invariant bilinear form induced from
Poincar\'e duality. To compute the image of monodromy, it is often
helpful to know this form explictly. After a review of the relevant
results of \cite{dw04} in Section 1, we give a formula for the
Poincar\'e duality pairing on $\W$ in Section 2. Finally, in Section 3
we illustrate our method in a very classical example: the
Picard--Euler system.


\section{Variation of parabolic cohomology revisited} \label{revise}

\subsection{}  \label{revise1}

Let $X$ be a compact Riemann surface of genus $0$ and $D\subset X$ a
subset of cardinality $r\geq 3$. We set $U:=X-D$. There exists a
homeomorphism $\kappa:X\iso\PP^1(\CC)$ between $X$ and the Riemann
sphere which maps the set $D$ to the real line $\PP^1(\RR)
\subset\PP^1(\CC)$. Such a homeomorphism is called a {\em marking} of
$(X,D)$.

Having chosen a marking $\kappa$, we may assume that $X=\PP^1(\CC)$
and $D\subset\PP^1(\RR)$. Choose a base point $x_0\in U$ lying in the
upper half plane. Write $D=\{x_1,\ldots,x_r\}$ with
$x_1<x_2<\ldots<x_r\leq\infty$. For $i=1,\ldots,r-1$ we let $\gamma_i$
denote the open interval $(x_i,x_{i+1})\subset U\cap\PP^1(\RR)$; for
$i=r$ we set $\gamma_0=\gamma_r:=(x_r,x_1)$ (which may include
$\infty$). For $i=1,\ldots,r$, we
let $\alpha_i\in\pi_1(U)$ be the element represented by a closed loop
based at $x_0$ which first intersects $\gamma_{i-1}$ and then
$\gamma_i$. We obtain the following well known presentation
\begin{equation} \label{pi1pres}
   \pi_1(U,x_0) \;=\; \gen{\;\alpha_1,\ldots,\alpha_r \,\mid\,
                                  \prod_i\alpha_i = 1\;},
\end{equation}
which only depends on the marking $\kappa$.

Let $R$ be a (commutative) ring. A {\em local system of $R$-modules}
on $U$ is a locally constant sheaf $\V$ on $U$ with values in the
category of free $R$-modules of finite rank. Such a local system
corresponds to a representation $\rho:\pi_1(U,x_0)\to\GL(V)$, where
$V:=\V_{x_0}$ is the stalk of $\V$ at $x_0$. For $i=1,\dots,r$, set
$g_i:=\rho(\alpha_i)\in\GL(V)$.  Then we have
\[
      \prod_{i=1}^r \; g_i \;=\; 1,
\]
and $\V$ can also be given by a tuple $\g=(g_1,\ldots,g_r)\in\GL(V)^r$
satisfying the above product-one-relation.

\begin{convention} \label{conv}
  Let $\alpha,\beta$ be two elements of $\pi_1(U,x_0)$, represented by
  closed path based at $x_0$. The composition $\alpha\beta$ is (the
  homotopy class of) the closed path obtained by first walking along
  $\alpha$ and then along $\beta$.  Moreover, we let $\GL(V)$ act on
  $V$ {\em from the right}.
\end{convention}

\subsection{}

Fix a local system of $R$-modules $\V$ on $U$ as above. Let $j:U\inj X$ denote
the inclusion. The {\em parabolic cohomology} of $\V$ is defined as the sheaf
cohomology of $j_*\V$, and is written as $ H^n\para(U,\V):=H^n(X,j_*\V)$. We
have natural morphisms $H_c^n(U,\V)\to H\para^n(U,\V)$ and $H\para^n(U,\V)\to
H^n(U,\V)$ ($H_c$ denotes cohomology with compact support).  Moreover, the group
$H^n(U,\V)$ is canonically isomorphic to the group cohomology
$H^n(\pi_1(U,x_0),V)$ and $H\para^1(U,\V)$ is the image of the cohomology with
compact support in $H^1(U,\V)$, see \cite{dw04}, Prop. 1.1.  Thus, there is a
natural inclusion
\[
      H^1\para(U,\V) \;\inj\; H^1(\pi_1(U,x_0),V).
\]

Let $\delta:\pi_1(U)\to V$ be a cocycle, i.e.\ we have
$\delta(\alpha\beta)=\delta(\alpha)\cdot\rho(\beta)+\delta(\beta)$ (see
Convention \ref{conv}). Set $v_i:=\delta(\alpha_i)$. It is clear that the tuple
$(v_i)$ is subject to the relation
\begin{equation} \label{vcocyclerel}
   v_1\cdot g_2\cdots g_r\,+\,v_2\cdot g_3\cdots g_r\,+\ldots+\,v_r \;=\; 0.
\end{equation} By definition, $\delta$ gives rise to 
an element in $H^1(\pi_1(U,x_0),V).$ We say that $\delta$ is a {\em parabolic}
cocycle if the class of $\delta$ in $H^1(\pi_1(U),V)$ lies in 
$H^1\para(U,\V)$. By \cite{dw04}, Lemma 1.2, the cocycle $\delta$ is parabolic
if and only if $v_i$ lies in the image of $g_i-1$, for all $i$.  Thus, the
association $\delta\mapsto(\delta(\alpha_1),\ldots, \delta(\alpha_r))$ yields an
isomorphism
\begin{equation} \label{Wgeq}
   H^1\para(U,\V) \;\cong\; W_\g:=H_\g/E_\g,
\end{equation}
where
\begin{equation} \label{Hgeq}
   H_{\g} \;:=\; \{\,(v_1,\ldots,v_r) \;\mid\;
     v_i\in{\rm Im}(g_i-1),\; \text{\rm relation \eqref{vcocyclerel} holds}\;\}
\end{equation}
and
\begin{equation} \label{Egeq}
   E_{\g} \;:=\; \{\, (\,v\cdot(g_1-1),\ldots,v\cdot(g_r-1)\,) \;\mid\;
                                    v\in V\;\}.
\end{equation}

\subsection{}  \label{locals3}

Let $S$ be a connected complex manifold, and $r\geq 3$. An {\em
$r$-configuration} over $S$ consists of a smooth and proper morphism $\pib:X\to
S$ of complex manifolds together with a smooth relative divisor $D\subset X$
such that the following holds. For all $s\in S$ the fiber $X_s:=\pib^{-1}(s)$ is
a compact Riemann surface of genus $0$. Moreover, the natural map $D\to S$ is an
unramified covering of degree $r$. Then for all $s\in S$ the divisor $D\cap X_s$
consists of $r$ pairwise distinct points $x_1,\ldots,x_r\in X_s$.

Let us fix an $r$-configuration $(X,D)$ over $S$.  We set $U:=X-D$ and denote by
$j:U\inj X$ the natural inclusion.  Also, we write $\pi:U\to S$ for the natural
projection.  Choose a base point $s_0\in S$ and set $X_0:=\pib^{-1}(s_0)$ and
$D_0:=X_0\cap D$. Set $U_0:=X_0-D_0=\pi^{-1}(s_0)$ and choose a base point
$x_0\in U_0$. The projection $\pi:U\to S$ is a topological fibration and yields
a short exact sequence
\begin{equation} \label{fibses}
  1 \;\To\; \pi_1(U_0,x_0) \;\To\; \pi_1(U,x_0) \;\To\; \pi_1(S,s_0) 
      \;\To\; 1.
\end{equation}

Let $\V_0$ be a local system of $R$-modules on $U_0$. A {\em variation} of
$\V_0$ over $S$ is a local system $\V$ of $R$-modules on $U$ whose restriction
to $U_0$ is identified with $\V_0$. The {\em parabolic cohomology} of a
variation $\V$ is the higher direct image sheaf
  \[
            \W \;:=\; R^1\pib_*(j_*\V).
  \]
By construction,  $\W$ is a local system with fibre
\[
     W \;:=\;  H^1\para(U_0,\V_0).
\] 
Thus $\W$ corresponds to a representation $\eta:\pi_1(S,s_0)\to\GL(W)$.  We call
$\rho$ the {\em monodromy representation} on the parabolic cohomology of $\V_0$
(with respect to the variation $\V$).

\subsection{}

Under a mild assumption, the monodromy representation $\eta$ has a very explicit
description in terms of the {\em Artin braid group}. We first have to introduce
some more notation. Define
\[
       \OO_{r-1} \;:=\; \{\, D'\subset \CC \;\mid\;
            |D'|=r-1\,\} \;=\; \{\,D\subset\PP^1(\CC) \;\mid\;
            |D|=r,\,\infty\in D\,\}.    
\]
The fundamental group $A_{r-1}:=\pi_1(\OO_{r-1},D_0)$ is the {\em Artin
  braid group} on $r-1$ strands. Let $\beta_1,\ldots,\beta_{r-2}$ be the
standard generators, see e.g.\ \cite{dw04}, \S 2.2. (The element $\beta_i$
switches the position of the two points $x_i$ and $x_{i+1}$; the point $x_i$
walks through the lower half plane and $x_{i+1}$ through the upper half plane.)
The generators $\beta_i$ satisfy the following well known relations:
\begin{equation} \label{braidrel}
  \beta_i\beta_{i+1}\beta_i \;=\; \beta_{i+1}\beta_i\beta_{i+1},
   \qquad \beta_i\beta_j \;=\; \beta_j\beta_i\;\;\;\text{\rm (for $|i-j|>1$).}
\end{equation}

Let $R$ be a commutative ring and $V$ a free $R$-module of finite
rank. Set
\[
    \E_r(V) \;:=\; \{\;\g=(g_1,\ldots,g_r) \mid 
                  g_i\in\GL(V),\;\; \prod_ig_i=1 \;\}.
\]
We define a right action of the Artin braid group $A_{r-1}$ on the set $\E_r(V)$ by
the following formula:
\begin{equation}  \label{locals5eq4}
   \g^{\beta_i} \;:=\; (g_1,\ldots,g_{i+1},g_{i+1}^{-1}g_ig_{i+1},\ldots,g_r).
\end{equation}
One easily checks that this definition is compatible with the relations
\eqref{braidrel}. For $\g\in\E_r(V)$, let $H_\g$ be as in \eqref{Hgeq}.  For all
$\beta \in A_{r-1}$, we define an $R$-linear isomorphism
\[
      \Phi(\g,\beta):\,H_{\g} \;\liso\; H_{\g^{\beta}},
\] 
as follows.  For the generators $\beta_i$ we set 
\begin{multline} \label{locals5eq6}
 \qquad (v_1,\ldots,v_r)^{\Phi(\g,\beta_i)} \\ 
    \;:=\;  (v_1,\ldots,v_{i+1},\,
    \underbrace{v_{i+1}(1-g_{i+1}^{-1}g_ig_{i+1})+v_ig_{i+1}}_{
        \text{\rm $(i+1)$th entry}},\,\ldots,v_r).
\end{multline} 
For an arbitrary word $\beta$ in the generators $\beta_i$, we define
 $\Phi(\g,\beta)$ using \eqref{locals5eq6} and the `cocycle rule'
\begin{equation} \label{locals5eq7}
    \Phi(\g,\beta)\cdot\Phi(\g^{\beta},\beta') \;=\; 
       \Phi(\g,\beta\beta').
\end{equation}
(Our convention is to let linear maps act from the right; therefore, the left
hand side of \eqref{locals5eq6} is the linear map obtained from first applying
$\Phi(\g,\beta)$ and then $\Phi(\g^\beta,\beta')$.)  It is easy to see that
$\Phi(\g,\beta)$ is well defined and respects the submodule $E_\g\subset H_\g$
defined by \eqref{Egeq}. Let
\[
       \bar{\Phi}(\g,\beta):\,W_{\g} \;\liso\; W_{\g^{\beta}}
\]
denote the induced map on the quotient $W_\g=H_{\g}/E_{\g}$. 

Given $\g\in\E_r(V)$ and $h\in\GL(V)$, we define the isomorphism
\[
   \Psi(\g,h):\; \left\{\;
     \begin{array}{ccc}
        H_{\g^h}  &  \;\liso\; & H_{\g} \\
        (v_1,\ldots,v_r) & \;\longmapsto\; & 
             (v_1\cdot h,\ldots,v_r\cdot h).
     \end{array}\right.,
\]
where $\g^h:=(h^{-1}g_1h,\ldots,h^{-1}g_r h)$. It is clear that $\Psi(\g,h)$
maps $E_{\g^h}$ to $E_\g$ and therefore induces an isomorphism
$\bar{\Psi}(\g,h):W_{\g^h}\iso W_\g$.

Note that the computation of the maps $\bar{\Phi}(\g,\beta)$ and
$\bar{\Psi}(\g,h)$ can easily be implemented on a computer.

\subsection{}  \label{mono3}

Let $S$ be a connected complex manifold, $s_0\in S$ a base point and $(X,D)$ an
$r$-configuration over $S$. As before we set $U:=X-D$, $D_0:=D\cap X_{s_0}$ and
$U_0:=U\cap X_{s_0}$. Let $\V_0$ be a local system of $R$-modules on $U_0$ and
$\V$ a variation of $\V_0$ over $S$. Let $\W$ be the parabolic cohomology of the
variation $\V$ and let $\eta:\pi_1(S,s_0)\to\GL(W)$ be the corresponding
monodromy representation. In order to describe $\eta$ explicitly, we find it
convenient to make the following assumption on $(X,D)$:

\begin{ass} \label{configass}
  \begin{enumerate}
  \item
    $X=\PP^1_S$ is the relative projective line over $S$.
  \item
    The divisor $D$ contains the section $\infty\times S\subset\PP^1_S$.
  \item There exists a point $s_0\in S$ such that $D_0:=D\cap\pib^{-1}(s_0)$ is
    contained in the real line $\PP^1(\RR)\subset\PP^1(\CC)=\pib^{-1}(s_0)$.
  \end{enumerate}
\end{ass}

In practise, this assumption is not a big restriction. See \cite{dw04} for a more
general setup. 

By Assumption \ref{configass}, we can consider $D_0$ as an element of
$\OO_{r-1}$.  Moreover, the divisor $D\subset\PP^1_S$ gives rise to an analytic
map $S\to\OO_{r-1}$ which sends $s_0\in S$ to $D_0\in\OO_{r-1}$. We let
$\varphi:\pi_1(S,s_0)\to A_{r-1}$ denote the induced group homomorphism and call
it the {\em braiding map} induced by $(X,D)$.

For $t\in\RR^+$ let $\Omega_t:=\{\,z\in\CC\,\mid\, |z|>t,\,
z\not\in(-\infty,0)\,\}$. Since $\Omega_t$ is contractible, the
fundamental group $\pi_1(U_0,\Omega_t)$ is well defined for $t\gg 0$
and independent of $t$, up to canonical isomorphism. We write
$\pi_1(U_0,\infty):=\pi_1(U_0,\Omega_t)$. We can define
$\pi_1(U,\infty)$ in a similar fashion, and obtain a short exact
sequence
\begin{equation} \label{fibses2}
  1 \;\to\; \pi_1(U_0,\infty) \;\To\; \pi_1(U,\infty) \;\To\; \pi_1(S,s_0) 
      \;\to\; 1.
\end{equation}
It is easy to see that the projection $\pi:U\to S$ has a continuous section
$\zeta:S\to U$ with the following property. For all $s\in S$ there exists $t\gg
0$ such that the region $\Omega_t$ is contained in the fibre
$U_s:=\pi^{-1}(s)\subset\PP^1(\CC)$ and such that $\zeta(s)\in\Omega_t$. The
section $\zeta$ induces a splitting of the sequence \eqref{fibses2}, which is
actually independent of $\zeta$. We will use this splitting to consider
$\pi_1(S,s_0)$ as a subgroup of $\pi_1(U,\infty)$. 

The variation $\V$ corresponds to a group homomorphism
$\rho:\pi_1(U,\infty)\to\GL(V)$, where $V$ is a free $R$-module. Let
$\rho_0$ denote the restriction of $\rho$ to $\pi_1(U_0,\infty)$ and
$\chi$ the restriction to $\pi_1(S,s_0)$. By Part (iii) of Assumption
\ref{configass} and the discussion in \S \ref{revise1} we have a
natural ordering $x_1<\ldots<x_r=\infty$ of the points in $D_0$, and a
natural choice of a presentation
$\pi_1(U_0,\infty)\cong\gen{\alpha_1,\ldots,\alpha_r\mid
\prod_i\alpha_i=1}$.  Therefore, the local system $\V_0$ corresponds
to a tuple $\g=(g_1,\ldots,g_r)\in\E_r(V)$, with
$g_i:=\rho_0(\alpha_i)$. One checks that the homomorphism
$\chi:\pi_1(S,s_0)\to\GL(V)$ satisfies the condition
\begin{equation} \label{chieq}
     \g^{\varphi(\gamma)} \;=\; \g^{\chi(\gamma)^{-1}},
\end{equation}
for all $\gamma\in\pi_1(S,s_0)$. Conversely, given $\g\in\E_r(V)$ and a
homomorphism $\chi:\pi_1(S,s_0)$ such that \eqref{chieq} holds then there exists
a unique variation $\V$ which induces the pair $(\g,\chi)$.

With these notations one has the following result (see 
\cite{dw04}, Thm. 2.5):

\begin{thm} \label{etathm} 
  Let $\W$ be the parabolic cohomology of $\V$ and
  $\eta:\pi_1(S,s_0)\to\GL(W_{\g})$ the corresponding monodromy
  representation. For all $\gamma\in\pi_1(S,s_0)$ we have
  \[
      \eta(\gamma)  \;=\;
           \bar{\Phi}(\g,\varphi(\gamma))\cdot\bar{\Psi}(\g,\chi(\gamma)).
   \]
\end{thm}

Thus, in order to compute the monodromy action on the
parabolic cohomology of a local system $\V_0$ corresponding to a tuple
$\g\in\E_r(V)$, we need to know the braiding map $\varphi:\pi_1(S,s_0)\to
A_{r-1}$ and the homomorphism $\chi:\pi_1(S,s_0)\to\GL(V)$. 

\begin{rem}
  Suppose that $R$ is a field and that the local system $\V_0$ is irreducible,
  i.e.\ the subgroup of $\GL(V)$ generated by the elements $g_i$ acts
  irreducibly on $V$. Then the homomorphism $\chi$ is determined, modulo the
  scalar action of $R^\times$ on $V$, by $\g$ and $\varphi$ (via
  \eqref{chieq}). It follows from Theorem \ref{etathm} that the projective
  representation $\pi_1(S,s_0)\to\PGL(V)$ associated to the monodromy
  representation $\eta$ is already determined by (and can be computed from) $\g$
  and the braiding map $\varphi$.
\end{rem}

\section{Poincar\'e duality}

Let $\V$ be a local system of $R$-modules on the punctured Riemann sphere
$U$. If $\V$ carries a non-degenerate symmetric (resp.\ alternating) form,
then Poincar\'e duality induces on the parabolic cohomology group
$H^1\para(U,\V)$ a non-degenerate alternating (resp.\ symmetric)
form. Similarly, if $R=\CC$ and $\V$ carries a Hermitian form, then we get a
Hermitian form on $H^1\para(U,\V)$.  In this section we derive an explicit
expression for this induced form.

\subsection{}
Let us briefly recall the definition of {\em singular (co)homology
with coefficients in a local system}. See e.g.\ \cite{Spanier93} for
more details. For $q\geq 0$ let $\Delta^q=|y_0,\ldots,y_q|$ denote the
standard $q$-simplex with vertices $y_0,\ldots,y_q$. We will sometimes
identify $\Delta^1$ with the closed unit interval $[0,1]$. Let $X$ be
a connected and locally contractible topological space and $\V$ a
local system of $R$-modules on $X$. For a continuous map $f:Y\to X$ we
denote by $\V_f$ the group of global sections of $f^*\V$.

In the following discussion, a {\em $q$-chain} will be a function $\varphi$
which assigns to each singular $q$-simplex $\sigma:\Delta^q\to X$ a section
$\varphi(\sigma)\in\V_{\sigma}$. Let $\Delta^q(X,\V)$ denote the set of all
$q$-chains, which is made into an $R$-module in the obvious way. A $q$-chain
$\varphi$ is said to have {\em compact support} if there exists a compact
subset $A\subset X$ such that $\varphi_\sigma=0$ whenever
$\supp(\sigma)\subset X-A$. The corresponding $R$-module is denoted by
$\Delta^q_c(X,\V)$. We define coboundary operators
$d:\Delta^q(X,\V)\to\Delta^{q+1}(X,\V)$ and
$d:\Delta^q_c(X,\V)\to\Delta^{q+1}_c(X,\V)$ through the formula
\[
    (d\, \varphi)(\sigma) \;:=\; \sum_{0\leq i\leq q}\,
       (-1)^i\cdot\overline{\varphi(\sigma^{(i)})}.
\]
Here $\sigma^{(i)}$ is the $i$th face of $\sigma$ (see \cite{Spanier})
and $\overline{\varphi(\sigma^{(i)})}$ denotes the unique extension of
$\varphi(\sigma^{(i)})$ to an element of $\V_\sigma$. It is proved in
\cite{Spanier93} that we have canonical isomorphisms
\begin{equation} \label{coheq1}
    H^n(X,\V) \;\cong\; H^n(\Delta^\bullet(X,\V),\,d), \qquad
    H^n_c(X,\V) \;\cong\; H^n(\Delta^\bullet_c(X,\V),\,d),
\end{equation}
i.e.\ singular cohomology agrees with sheaf cohomology. Let $x_0\in X$ be a base
point and $V$ the fibre of $\V$ at $x_0$. Then we also have an isomorphism
\begin{equation} \label{coheq2}
   H^1(X,\V) \;\cong\; H^1(\pi_1(X,x_0),V).
\end{equation}
Let $\varphi$ be a $1$-chain with $d\varphi=0$. Let $\alpha:[0,1]\to X$ be a
closed path with base point $x_0$.  By definition, $\varphi(\alpha)$ is a global
section of $\alpha^*\V$. Then $\alpha\mapsto\delta(\alpha):=\varphi(\alpha)(1)$
defines a cocycle $\delta:\pi_1(X,x_0)\to V$, and this cocycle represents the
image of $\varphi$ in $H^1(X,\V)$.

A $q$-chain $\varphi$ is called {\em finite} if $\varphi(\sigma)=0$ for all
but finitely many simplexes $\sigma$. It is called {\em locally finite} if
every point in $X$ has a neighborhood $U\subset X$ such that
$\varphi(\sigma)=0$ for all but finitely many simplexes $\sigma$ contained in
$U$. We denote by $\Delta_q(X,\V)$ (resp.\ by $\Delta_q^{lf}(X,\V)$) the
$R$-module of all finite (resp.\ locally finite) $q$-chains. For a fixed
$q$-simplex $\sigma$ and a section $v\in\V_\sigma$, the symbol
$v\otimes\sigma$ will denote the $q$-chain which assigns $v$ to $\sigma$ and
$0$ to all $\sigma'\not=\sigma$.  Obviously, every finite (resp.\ locally
finite) $q$-chain can be written as a finite (resp.\ possibly infinite) sum
$\sum_\mu\,v_\mu\otimes\sigma_\mu$. We define boundary operators
$\partial:\Delta_q(X,\V)\to\Delta_{q-1}(X,\V)$ and
$\partial:\Delta_q^{lf}(X,\V)\to\Delta_{q-1}^{lf}(X,\V)$ through the formula
\[
   \partial(v\otimes\sigma) \;:=\; \sum_{0\leq i\leq q}\,
       (-1)^i\cdot v|_{\sigma^{(i)}}\otimes\sigma^{(i)}.
\]
We define homology (resp.\ locally finite homology) with coefficients in $\V$
as follows:
\[
    H_q(X,\V) \;:=\; H_q(\Delta_\bullet(X,\V)), \qquad
    H_q^{lf}(X,\V) \;:=\; H_q(\Delta_\bullet^{lf}(X,\V)).
\]

\subsection{}  \label{poincare2}        

Let $X:=\PP^1(\CC)$ be the Riemann sphere and
$D=\{x_1,\ldots,x_r\}\subset\PP^1(\RR)$ a subset of $r\geq 3$ points
lying on the real line, with $x_1<\ldots<x_r\leq\infty$. Let $\V$ be a
local system of $R$-modules on $U=X-D$. Choose a base point $x_0$
lying in the upper half plane. Then $\V$ corresponds to a tuple
$\g=(g_1,\ldots,g_r)$ in $\GL(V)$ with $\prod_ig_i=1$, where
$V:=\V_{x_0}$. See \S \ref{revise1}. Let $\V^*:=\underline{\rm
Hom}(\V,R)$ denote the local system dual to $\V$. It corresponds to
the tuple $\g^*=(g_1^*,\ldots,g_r^*)$ in $\GL(V^*)$, where $V^*$ is
the dual of $V$ and for each $g\in\GL(V)$ we let $g^*\in\GL(V^*)$ be
the unique element such that
\[
     \gen{w\cdot g^*,v\cdot g} \;=\; \gen{w,v}
\]
for all $w\in V^*$ and $v\in V$. Note that $V^{**}=V$ because $V$ is free of
finite rank over $R$. 

Let $\varphi$ be a $1$-chain with compact support and with coefficients in
$\V^*$. Let $a=\sum_\mu v_\mu\otimes\alpha_\mu$ be a locally finite $1$-chain
with coefficients in $\V$. By abuse of notation, we will also write $\varphi$
(resp.\ $a$) for its class in $H^1_c(U,\V^*)$ (resp.\ in $H_1^{lf}(U,\V)$).
The {\em cap product}
\[
     \varphi \cap a \;:=\; \sum_\mu\,\gen{\varphi(\alpha_\mu),v_\mu}
\]
induces a bilinear pairing
\begin{equation} \label{poincareeq2}
     \cap\,:\; H^1_c(U,\V^*) \;\otimes\; H_1^{lf}(U,\V) \;\To\; R.
\end{equation}
It is easy to see from the definition that $H_0^{lf}(U,\V)=0$. Therefore, it
follows from the Universal Coefficient Theorem for cohomology (see e.g.\ 
\cite{Spanier}, Thm.\ 5.5.3) that the pairing \eqref{poincareeq2} is
nonsingular on the left, i.e.\ identifies $H^1_c(U,\V^*)$ with ${\rm
  Hom}\,(H_1^{lf}(U,\V),R)$. The cap product also induces a pairing
\begin{equation} \label{poincareeq3}
     \cap\,:\; H^1(U,\V^*) \;\otimes\; H_1(U,\V) \;\To\; R.
\end{equation}
(This last pairing may not be non-singular on the left. The
reason is that
\[
         H_0(U,\V) \;\cong\; V/\gen{\,{\rm Im}(g_i-1)\mid i=1,\ldots,r\,} 
\]
may not be a free $R$-module, and so ${\rm Ext}^1(H_0(U,\V),R)$ may be
nontrivial.) Let $f^1:H^1_c(U,\V^*)\to H^1(U,\V^*)$ and $f_1:H_1(U,\V)\to
H_1^{lf}(U,\V)$ denote the canonical maps. Going back to the definition, one
can easily verify the rule
\begin{equation} \label{poincareeq4}
   f^1(\varphi)\cap a \;=\; \varphi\cap f_1(a).
\end{equation}

Let $\varphi\in H^1_c(U,\V^*)$ and $\psi\in H^1(U,\V)$. The {\em cup
product} $\varphi\cup\psi$ is defined as an element of $H^2_c(U,R)$,
see \cite{Steenrod43} or \cite{Spanier93}. The standard orientation of
$U$ yields an isomorphism $H^2_c(U,R)\cong R$. Using this isomorphism,
we shall view the cup product as a bilinear pairing
\[
    \cup\,:\;H^1_c(U,\V^*)\otimes H^1(U,\V) \;\To\; R.
\]
Similarly, one can define the cup product $\varphi\cup\psi$, where $\varphi\in
H^1(U,\V^*)$ and $\psi\in H^1_c(U,\V)$. Given $\varphi\in H^1_c(U,\V^*)$ and
$\psi\in H^1_c(U,\V)$, one checks that 
\begin{equation} \label{poincareeq5}
  f^1(\varphi)\cup\psi \;=\; \varphi\cup f^1(\psi).
\end{equation}

\begin{prop}[Poincar\'e duality] \label{poincareprop}
  There exist unique isomorphisms of $R$-modules
  \[
     p:\, H_1(U,\V) \liso H^1_c(U,\V), \qquad
     p:\, H_1^{lf}(U,\V) \liso H^1(U,\V)
  \]
  such that the following holds. If $\varphi\in H^1_c(U,\V^*)$ and $a\in
  H_1^{lf}(U,\V)$ or if $\varphi\in H^1(U,\V^*)$ and $a\in
  H_1(U,\V)$ then we have
  \[
       \varphi \cap a \;=\; \varphi \cup p(a).
  \]
  These isomorphisms are compatible with the canonical maps $f_1$ and $f^1$,
  i.e.\ we have $p\circ f_1=f^1\circ p$. 
\end{prop}

\proof
See \cite{Steenrod43} or \cite{Spanier93}.
\Endproof

\begin{cor} \label{poincarecor}
  The cup product induces a non-degenerate bilinear pairing
  \[
     \cup\,:\; H^1\para(U,\V^*) \otimes H^1\para(U,\V)\;\To\; R.
  \]
\end{cor}

\proof Let $\varphi\in H^1\para(U,\V^*)$ and $\psi\in H^1\para(U,\V)$. Choose
$\varphi'\in H^1_c(U,\V^*)$ and $\psi'\in H^1_c(U,\V)$ with
$\varphi=f^1(\varphi')$ and $\psi=f^1(\psi')$.  By \eqref{poincareeq5} we have
$\varphi'\cup\psi=\varphi\cup\psi'$. Therefore, the expression
$\varphi\cup\psi:=\varphi'\cup\psi$ does not depend on the choice of the lift
$\varphi'$ and defines a bilinear pairing between $H^1\para(U,\V^*)$ and
$H^1\para(U,\V)$. By Proposition \ref{poincareprop} and since the cap product
\eqref{poincareeq2} is non-degenerate on the left, this pairing is also
non-degenerate on the left. But the cup product is alternating (i.e.\ we have
$\varphi\cup\psi=-\psi\cup\varphi$, where the right hand side is defined using
the identification $\V^{**}=\V$), so our pairing is also non-degenerate on the
right.    \Endproof

For $a\in H_1^{lf}(U,\V^*)$ and $b\in H_1(U,\V)$, the expression
\[
       (a,b) \;:=\; p(a)\cup p(b)
\]
defines another bilinear pairing $H_1^{lf}(U,\V^*)\otimes H_1(U,\V)\to
R$. It is shown in \cite{Steenrod43} that this pairing can be computed
as an `intersection product of loaded cycles', generalizing the usual
intersection product for constant coefficients, as follows. We may
assume that $a$ is represented by a locally finite chain $\sum_\mu
v^*_\mu\otimes\alpha_\mu$ and that $b$ is represented by a finite
chain $\sum_\nu v_\nu\otimes\beta_\nu$ such that for all $\mu,\nu$ the
$1$-simplexes $\alpha_\mu$ and $\beta_\nu$ are smooth and intersect
each other transversally, in at most finitely many points. Suppose $x$
is a point where $\alpha_\mu$ intersects $\beta_\nu$. Then there
exists $t_0\in [0,1]$ such that $x=\alpha(t_0)=\beta(t_0)$ and
$(\frac{\partial\alpha}{\partial
t}|_{t_0},\frac{\partial\beta}{\partial t}|_{t_0})$ is a basis of the
tangent space of $U$ at $x$. We set $\imath(\alpha,\beta,x):=1$
(resp.\ $\imath(\alpha,\beta,x):=-1$) if this basis is positively
(resp.\ negatively) oriented. Furthermore, we let $\alpha_{\mu,x}$
(resp. $\beta_{\nu,x}$) be the restriction of $\alpha$ (resp. of $\beta$)
to the interval $[0,t_0]$. Then we have
\begin{equation} \label{poincareeq6}
  (a,b) \;=\; \sum_{\mu,\nu,x}\; \imath(\alpha_\mu,\beta_\nu,x)\cdot
                 \gen{\,(v^*)^{\alpha_{\mu,x}},v^{\beta_{\nu,x}}\,}.
\end{equation}

\subsection{}

Let $\V\otimes\V\to\underline{R}$ be a non-degenerate symmetric (resp.\ 
alternating) bilinear form, corresponding to an injective homomorphism
$\kappa:\V\inj\V^*$ with $\kappa^*=\kappa$ (resp.\ $\kappa^*=-\kappa$). We
denote the induced map $H^1\para(U,\V)\to H^1\para(U,\V^*)$ by $\kappa$ as well.
Then
\[
        \gen{\varphi,\psi} \;:=\; \kappa(\varphi)\cup\psi
\]
defines a non-degenerate alternating (resp.\ symmetric)
form on $H^1\para(U,\V)$.

Similarly, suppose that $R=\CC$ and let $\V$ be equipped with a non-degenerate
Hermitian form, corresponding to an isomorphism $\kappa:\VB\iso\V^*$. Then the
pairing
\begin{equation} \label{hermform}
   (\varphi,\psi) \;:=\; -i\cdot(\kappa(\bar{\varphi})\cup\psi)
\end{equation}
is a nondegenerate Hermitian form on $H^1_p(U,\V)$ (we identify $H^1_p(U,\VB)$
with the complex conjugate of the vector space $H^1_p(U,\V)$ in the obvious
way). 

Suppose that the Hermitian form on $\V$ is positive definite. Then we
can express the signature of the form \eqref{hermform} in terms of the
tuple $\g$, as follows. For $i=1,\ldots,r$, let
\[
   g_i \;\sim\;\; \begin{pmatrix} 
      \alpha_{i,1} & & \\ & \ddots & \\ & & \alpha_{i,n} 
                  \end{pmatrix}
\]
be a diagonalization of $g_i\in\GL(V)$. Since the $g_i$ are Hermitian,
the eigenvalues $\alpha_{i,j}$ have absolute value one and can be
uniquely written in the form $\alpha_{i,j}=\exp(2\pi i \mu_{i,j})$,
with $0\leq \mu_{i,j}<1$. Set $\bar{\mu}_{i,j}:=1-\mu_{i,j}$ if
$\mu_{i,j}>0$ and $\bar{\mu}_{i,j}:=0$ otherwise. 

\begin{thm}  \label{signaturethm}
  Suppose that $\V$ is equipped with a positive definite Hermitian
  form and that $H^0(U,\V)=0$. Then the Hermitian form
  \eqref{hermform} on $H^1\para(U,\V)$ has signature
  \[
     \big(\,(\sum_{i,j}\mu_{i,j})-\dim_\CC V,\, 
            (\sum_{i,j}\bar{\mu}_{i,j})-\dim_\CC V\,\big).
  \]
\end{thm}

If $\dim_\CC V=1$, this formula is proved in \cite{DeligneMostow}, \S
2. With some extra work, the proof can be generalized to the case of
arbitrary dimension. See forthcoming work of the authors.

\subsection{}

We are interested in an explicit expression for the pairing of
Corollary \ref{poincarecor}. We use the notation introduced at the
beginning of \S \ref{poincare2}, with the following modification. By
$\gamma_i$ we now denote a homeomorphism between the open unit
interval $(0,1)$ and the open interval $(x_i,x_{i+1})$. We assume that
$\gamma_i$ extends to a path $\bar{\gamma}_i:[0,1]\to\PP^1(\RR)$ from
$x_i$ to $x_{i+1}$. We denote by $U^+\subset\PP^1(\CC)$ (resp. $U^-$)
the upper (resp.\ the lower) half plane and by $\bar{U}^+$ (resp.\
$\bar{U}^-$) its closure inside $U=\PP^1(\CC)-\{x_1,\ldots,x_r\}$.
Since $\bar{U}^+$ is simply connected and contains the base point
$x_0$, an element of $V$ extends uniquely to a section of $\V$ over
$\bar{U}^+$. We may therefore identify $V$ with $\V(\bar{U}^+)$ and
with the stalk of $\V$ at any point $x\in\bar{U}^+$.

Choose a sequence of numbers $\epsilon_n$, $n\in\ZZ$, with
$0<e_n<e_{n+1}<1$ such that $\epsilon_n\to 0$ for $n\to-\infty$ and
$\epsilon_n\to 1$ for $n\to\infty$.  Let $\gamma_i^{(n)}:[0,1]\to U$
be the path
$\gamma_i^{(n)}(t):=\gamma_i(\epsilon_nt+\epsilon_{n-1}(1-t))$.
Let $w_1,\ldots,w_r\in V$. Since
$\supp(\gamma_i)\subset\bar{U}^+$, it makes sense to define
\[
     w_i\otimes\gamma_i \;:=\; \sum_n\,w_i\otimes\gamma_i^{(n)}.
\]
This is a locally finite $1$-chain. Set 
\[
     c \;:=\; \sum_{i=1}^r\, w_i\otimes\gamma_i.
\]
Note that $\partial(c)=0$, so $c$ represents a class in $H_1^{lf}(U,\V)$. 

\begin{lem} \label{poincarelem}
\begin{enumerate}
\item
  The image of $c$ under the Poincar\'e isomorphism $H_1^{lf}(U,\V)\cong
  H^1(U,\V)$ is represented by the unique cocycle $\delta:\pi_1(U,x_0)\to V$
  with 
  \[
       \delta(\alpha_i) \;=\; w_i - w_{i-1}\cdot g_i.
  \]
\item
  The cocycle $\delta$ in (i) is parabolic if and only if there exist
  elements $u_i\in V$ with $w_i-w_{i-1}=u_i\cdot(g_i-1)$, for all $i$.
\end{enumerate}
\end{lem}

\proof For a path $\alpha:[0,1]\to U$ in $U$, consider the following 
conditions:
\begin{itemize}
\item[(a)]
  The support of $\alpha$ is contained either in $U^+$ or in $U^-$.
\item[(b)]
  We have $\alpha(0)\in U^+$, $\alpha(1)\in U^-$ and $\alpha$ intersects
  $\gamma_i$ transversally in a unique point.   
\item[(c)]
  We have $\alpha(0)\in U^-$, $\alpha(1)\in U^+$ and $\alpha$ intersects
  $\gamma_i$ transversally in a unique point.   
\end{itemize}
In Case (b) (resp.\ in Case (c)) we identify $\V_\alpha$ with $V$ via
the stalk $\V_{\alpha(0)}$ (resp.\ via $\V_{\alpha(1)}$.  Let
$\varphi\in C^1(U,\V)$ be the unique cocycle such that
\[
 \varphi(\alpha) \;=\;\; 
 \begin{cases}
   \;\; 0,   & \;\text{\rm if $\alpha$ is as in Case (a)} \\
   \;\; -w_i, & \;\text{\rm if $\alpha$ is as in Case (b)} \\
   \;\; w_i^{\alpha^{-1}},&\;\text{\rm if $\alpha$ is as in Case (c).}
 \end{cases}   
\]
(To show the existence and uniqueness of $\varphi$, choose a triangulation of
$U$ in which all edges satisfy Condition (a), (b) or (c). Then use simplicial
approximation.) We claim that $\varphi$ represents the image of the cycle $c$
under the Poincar\'e isomorphism. Indeed, this follows from the definition of
the Poincar\'e isomorphism, as it is given in \cite{Steenrod43}. Write
$\alpha_i=\alpha_i'\alpha_i''$, with $\alpha_i'(1)=\alpha_i''(0)\in U^-$. 
Using the fact that $\varphi$ is a cocycle we get
\[
    \varphi(\alpha_i) \;=\; 
      \varphi(\alpha_i') + \varphi(\alpha_i'')^{{\alpha_i'}^{-1}} \;=\;
           -w_{i-1} + w_i\cdot g_i^{-1}.
\]
Therefore we have $\delta(\alpha_i)=\varphi(\alpha_i)\cdot g_i=
w_i-w_{i-1}\cdot g_i$.  See Figure \ref{poincarefig1}. This proves (i).

\begin{figure}
\begin{center}

\setlength{\unitlength}{.27mm}
\begin{picture}(250,160)(-20,30)
\put(40,60){\circle*{4}}
\put(110,60){\circle*{4}}
\put(180,60){\circle*{4}}
\put(90,130){\circle*{4}}
\put(99,128){$\scriptstyle x_0$}
\put(40,60){\line(1,0){140}}
\put(75,60){\vector(1,0){2}}
\put(145,60){\vector(1,0){2}}
\dashline{5}(-20,60)(40,60)
\dashline{5}(180,60)(240,60)
\curve(90,130,95,50,120,35,125,70,90,130)
\put(90,80){\vector(0,-1){1}}
\put(121,80){\vector(-1,2){1}}
\put(37,48){$\scriptstyle x_{i-1}$}
\put(107,48){$\scriptstyle x_i$}
\put(173,48){$\scriptstyle x_{i+1}$}
\put(64,67){$\scriptstyle \gamma_{i-1}$}
\put(145,67){$\scriptstyle \gamma_i$}
\put(124,93){$\scriptstyle \alpha_i$}
\put(210,90){$U^+$}
\put(210,20){$U^-$}
\end{picture}

\end{center}
\caption{\label{poincarefig1}}
\end{figure}

By Section 1.1, the cocycle $\delta$ is parabolic if and only if
$v_i$ lies in the image of $g_i-1$. So (ii) follows from (i) by a simple
manipulation. \Endproof

\begin{thm} \label{poincarethm}
  Let $\varphi\in H^1\para(U,\V^*)$ and $\psi\in H^1\para(U,\V)$, represented
  by cocycles $\delta^*:\pi_1(U,x_0)\to V^*$ and $\delta:\pi_1(U,x_0)\to
  V$. Set $v_i:=\delta(\alpha_i)$ and $v_i^*=\delta^*(\alpha_i)$. 
  If we choose $v_i'\in V$ such that $v_i'\cdot(g_i-1)=v_i$ (see 
Lemma \ref{poincarelem}), then we have
  \[
     \varphi\cup\psi \;=\; \sum_{i=1}^r\;
       (\;\gen{v_i^*,v_i'}\;+\;\sum_{j=1}^{i-1}
          \gen{v_j^*g_{j+1}^*\cdots g_{i-1}^*(g_i^*-1),v_i'}\;).
  \]
\end{thm}

\proof
  Let $w_1:=v_1$, $w_1^*:=v_1^*$ and 
  \[
       w_i \;:=\; v_i + w_{i-1}\cdot g_i, \qquad
       w_i^* \;:=\; v_i^* + w_{i-1}^*\cdot g_i^*
  \]
  for $i=2,\ldots,r$. By Lemma \ref{poincarelem},
 we can choose $u_i\in V$ with
  $w_i-w_{i-1}=u_i\cdot(g_i-1)$, for $i=1,\ldots,r$. The
claim will follow from the following formula: 
  \begin{equation}\label{www}
     \varphi\cup\psi \;=\; \sum_{i=1}^r\; \gen{w_i^*-w_{i-1}^*,u_i-w_{i-1}}.
  \end{equation}
To prove Equation \eqref{www},
suppose $\delta$ is parabolic, and choose $u_i\in V$ such that
$w_i-w_{i-1}=u_i\cdot(g_i-1)$. Let $D_i\subset X$ be a closed disk
containing $x_i$ but none of the other points $x_j$, $j\not= i$. We may assume
that the boundary of $D_i$ intersects $\gamma_{i-1}$ in the point
$\gamma_{i-1}^{(0)}(1)$ but nowhere else, and that $D_i$ intersects $\gamma_i$
in the point $\gamma_i^{(0)}(0)$ but nowhere else. Set
$D_i^+:=D_i\cap\bar{U}^+$ and $D_i^-:=D_i\cap\bar{U}^-$. Let
$u_i^+:=u_i-w_{i-1}$, considered as a section of $\V$ over $D_i^+$ via
extension over the whole upper half plane $U^+$. It makes sense to define the
locally finite chain
\[
      u_i^+\otimes D_i^+ \;:=\; \sum_\sigma\, u_i^+\otimes\sigma,
\]
where $\sigma$ runs over all $2$-simplexes of a triangulation of $D_i^+$.
(Note that $x_i\not\in D_i^+$, so this triangulation cannot be finite.)
Similarly, let $u_i^-\in\V_{D_i^-}$ denote the section of $\V$ over $D_i^-$
obtained from $u_i\in V$ by continuation along a path which enters $U^-$
from $U^+$ by crossing the path $\gamma_{i-1}$; define $u_i^-\otimes
D_i^-$ as before. Let 
\[
    c' \;:=\; c + \partial\,(u_i^+\otimes D_i^+ + u_i^-\otimes D_i^-).
\]
It is easy to check that $c'$ is homologous to the cocycle
\[
  c'' \;:=\; \sum_i\;\big(\,w_i\otimes\gamma_i^{(0)} + 
      u_i^+\otimes\beta_i^+ + u_i^-\otimes\beta_i^-\,\big),
\]
where $\beta_i^+$ (resp.\ $\beta_i^-$) is the path from $\gamma_i^{(0)}(0)$ to
$\gamma_{i-1}^{(0)}(1)$ (resp.\ from $\gamma_{i-1}^{(0)}(1)$ to
$\gamma_i^{(0)}(0)$) running along the upper (resp.\ lower) part of the
boundary of $D_i$. See Figure \ref{poincarefig2}. Note that $c''$ is
finite and that, by construction, the image of $c''$ under the canonical map
$f_1:H_1(U,\V)\to H_1^{lf}(U,\V)$ is equal to the class of $c$. Let $\psi'\in
H^1_c(U,\V)$ denote the image of $c''$ under the Poincar\'e isomorphism
$H_1(U,\V)\cong H^1_c(U,\V)$. The last statement of Proposition
\ref{poincareprop} shows that $\psi'$ is a lift of $\psi\in H^1\para(U,\V)$.

\begin{figure}
\begin{center}

\setlength{\unitlength}{.35mm}
\begin{picture}(300,100)(0,10)
\put(50,50){\circle*{4}}
\put(150,50){\circle*{4}}
\put(250,50){\circle*{4}}
\put(50,50){\circle{40}}
\put(150,50){\circle{40}}
\put(250,50){\circle{40}}
\drawline(0,50)(30,50)
\drawline(70,50)(130,50)
\drawline(170,50)(230,50)
\drawline(270,50)(300,50)
\dottedline{2}(30,50)(70,50)
\dottedline{2}(130,50)(170,50)
\dottedline{2}(230,50)(270,50)
\curve(50,50,100,65,150,50)
\put(100,65){\vector(1,0){0}}
\put(95,72){$\scriptstyle \gamma_{i-1}'$}
\curve(150,50,200,65,250,50)
\put(200,65){\vector(1,0){0}}
\put(195,72){$\scriptstyle \gamma_i'$}
\put(45,39){$\scriptstyle x_{i-1}$}
\put(145,39){$\scriptstyle x_i$}
\put(245,39){$\scriptstyle x_{i+1}$}
\put(100,50){\vector(1,0){1}}
\put(93,35){$\scriptstyle \gamma_{i-1}^{(0)}$}
\put(200,50){\vector(1,0){1}}
\put(193,35){$\scriptstyle \gamma_i^{(0)}$}
\put(49,70){\vector(-1,0){0}}
\put(149,70){\vector(-1,0){0}}
\put(145,77){$\scriptstyle \beta_i^+$}
\put(249,70){\vector(-1,0){0}}
\put(51,30){\vector(1,0){0}}
\put(151,30){\vector(1,0){0}}
\put(145,18){$\scriptstyle \beta_i^-$}
\put(251,30){\vector(1,0){0}}
\end{picture}

\end{center}
\caption{\label{poincarefig2}}
\end{figure}

Let $c^*:=\sum_i w_i^*\otimes\gamma_i\in C_1(U,\V^*)$. By (i) and the choice
of $w_i^*$, the image of $c^*$ under the Poincar\'e isomorphism
$H_1^{lf}(U,\V^*)\cong H^1(U,\V^*)$ is equal to $\varphi$. By definition, we
have $\varphi\cup\psi=(c^*,c'')$. To compute this intersection number, we have
to replace $c^*$ by a homologous cycle which intersects the support of $c''$
at most transversally. For instance, we can deform the open paths $\gamma_i$
into open paths $\gamma_i'$ which lie entirely in the upper half plane. See
Figure \ref{poincarefig2}. It follows from \eqref{poincareeq6} that
\[
    (c^*,c'') \;=\; \sum_i\;\;
       \gen{w_{i-1}^*,u_i^+}-\gen{w_i^*,u_i^+} \;=\;
       \sum_i\;\; \gen{w_i^*-w_{i-1}^*,u_i-w_{i-1}}.
\]
This finishes the proof of \eqref{www}. The formula in (iv) follows from
\eqref{www} from a
straightforward computation, expressing $w_i$ and $u_i$ in terms of 
$v_i$ and $v_i'.$
\Endproof

\begin{rem}
  In the somewhat different setup, a similar formula as in Theorem \ref{poincarethm} 
  can be found in \cite{Voelklein01}, \S 1.2.3. 
\end{rem}


\section{The monodromy of the Picard--Euler system} \label{picard}

Let 
\[
    S \;:=\; \{\,(s,t)\in\CC^2 \;\mid\; s,t\not=0,1,\; s\not=t\;\},
\]
and let $X:=\PP^1_S$ denote the relative projective line over $S$. 
The equation 
\begin{equation} \label{picardeq}
       y^3 \;=\; x(x-1)(x-s)(x-t)
\end{equation}
defines a finite Galois cover $f:Y\to X$ of smooth projective curves
over $S$, tamely ramified along the divisor
$D:=\{0,1,s,t,\infty\}\subset X$. The curve $Y$ is called the {\em
Picard curve}. Let $G$ denote the Galois group of $f$, which is cyclic
of order $3$. The equation $\sigma^*y=\chi(\sigma)\cdot y$ for
$\sigma\in G$ defines an injective character
$\chi:G\inj\CC^\times$. As we will see below, the $\chi$-eigenspace of
the cohomology of $Y$ gives rise to a local system on $S$ whose
associated system of differential equations is known as the {\em
Picard--Euler system}.

We fix a generator $\sigma$ of $G$ and set $\omega:=\chi(\sigma)$. Let
$K:=\QQ(\omega)$ be the quadratic extension of $\QQ$ generated by
$\omega$ and $\OO_K=\ZZ[\omega]$ its ring of integers. The family of
$G$-covers $f:Y\to X$ together with the character $\chi$ of $G$
corresponds to a local system of $\OO_K$-modules on $U:=X-D$. Set
$s_0:=(2,3)\in S$ and let $\V_0$ denote the restriction of $\V$ to the
fibre $U_0=\AA^1_\CC-\{0,1,2,3\}$ of $U\to S$ over $s_0$.  We consider
$\V$ as a variation of $\V_0$ over $S$. Let $\W$ denote the parabolic
cohomology of this variation; it is a local system of $\OO_K$-modules
of rank three, see \cite{dw04}, Rem. 1.4.  Let $\chi':G\inj\CC^\times$
denote the conjugate character to $\chi$ and $\W'$ the parabolic
cohomology of the variation of local systems $\V'$ corresponding to
the $G$-cover $f$ and the character $\chi'$. We write $\W_\CC$ for the
local system of $\CC$-vectorspaces $\W\otimes\CC$. The maps
$\pi_Y:Y\to S$ and $\pi_X:X\to S$ denote the natural projections.

\begin{prop} \label{picardprop1}
  We have a canonical isomorphism of local systems
  \[
       R^1\pi_{Y,*}\underline{\CC} \;\cong\; 
            \W_\CC \,\oplus\, \W'_\CC.
  \]
  This isomorphism identifies the fibres of $\W_\CC$ with the
  $\chi$-eigenspace of the singular cohomology of the Picard curves
  of the family $f$.  
\end{prop}

\proof
The group $G$ has a natural left action on the sheaf $f_*\underline{\CC}$.
We have a canonical isomorphism of sheaves on $X$
\[
       f_*\underline{\CC} \;\cong\; \underline{\CC} \,\oplus\,
        j_*\V_\CC \,\oplus\, j_*\V',
\]
which identifies $j_*\V_\CC$, fibre by fibre, with the
$\chi$-eigenspace of $f_*\underline{\CC}$. Now the Leray spectral
sequence for the composition $\pi_Y=\pi_X\circ f$ gives
isomorphisms of sheaves on $S$
\[
   R^1\pi_{Y,*}\underline{\CC} \;\cong\;
     R^1\pi_{X,*}(f_*\underline{\CC}) \;\cong\;
       \W_\CC \,\oplus\, \W'_\CC.
\]
Note that $R^1\pi_{X,*}\underline{\CC}=0$ because the genus of $X$ is zero.
Since the formation of $R^1\pi_{Y,*}$ commutes with the $G$-action, the
proposition follows.  \Endproof


The comparison theorem between singular and deRham cohomology
identifies $R^1\pi_{Y,*}\underline{\CC}$ with the local system of
horizontal sections of the relative deRham cohomology module
$R^1\dR\pi_{Y,*}\OO_Y$, with respect to the Gauss-Manin connection. The
$\chi$-eigenspace of $R^1\dR\pi_{Y,*}\OO_Y$ gives rise to a Fuchsian
system known as the Picard--Euler system. In more classical terms, the
Picard--Euler system is a set of three explicit partial differential
equations in $s$ and $t$ of which the period integrals
\[
      I(s,t;a,b) \;:=\; 
        \int_a^b\,\frac{{\rm d}\,x}{\sqrt[3]{x(x-1)(x-s)(x-t)}}
\]
(with $a,b\in\{0,1,s,t,\infty\}$) are a solution. See \cite{Picard83},
\cite{HolzapfelEuler}, \cite{HolzapfelBall}. It follows from
Proposition \ref{picardprop1} that the monodromy of the Picard--Euler
system can be identified with the representation
$\eta:\pi_1(S)\to\GL_3(\OO_K)$ corresponding to the local system $\W$. 

\begin{thm}[Picard] \label{picardthm}
  For suitable generators $\gamma_1,\ldots,\gamma_5$ of the fundamental
  group $\pi_1(S)$, the matrices
  $\eta(\gamma_1),\ldots,\eta(\gamma_5)$ are equal to 
  \begin{gather*}
     \begin{pmatrix}
        \omega^2        & \;0\;          & 1-\omega            \\
        \omega-\omega^2 & 1              & \omega^2 -1         \\
        0               & 0              & 1 
     \end{pmatrix},\;
     \begin{pmatrix}
        \omega^2        & \;0\;          & 1-\omega^2          \\
        1-\omega^2      & 1              & \omega^2-1          \\
        0               & 0              & 1 
     \end{pmatrix},\;
     \begin{pmatrix}
        \;\;1\;\;       & 0              & 0                    \\
        0               & \omega         & \omega^2-1           \\
        0               & \omega^2-1     & -2\omega   
     \end{pmatrix},\\
     \begin{pmatrix}
        \;\;\omega^2\;  & \;\;0\;\;      & \;\;0\;\;            \\
        0               & 1              & 0                    \\
        0               & 0              & 1 
     \end{pmatrix},\;
     \begin{pmatrix}
        \omega^2        & \omega-\omega^2& \;0\;\\
        0               & 1              & 0    \\
        1-\omega        & \omega^2-1     & 1 
     \end{pmatrix}.
  \end{gather*} The invariant Hermitian form (induced by Poincar\'e 
duality, see Corollary \ref{poincarecor}) is given by the matrix
 \begin{gather*}
     \begin{pmatrix}
             a  & \;0\;          &    0        \\
        0 & 0              & a         \\
        0               & a           & 0 
     \end{pmatrix},\;
  \end{gather*}
where $a=\frac{i}{3}(\omega^2-\omega).$

\end{thm}     

\begin{figure}
\begin{center}

\setlength{\unitlength}{0.0011in}
\begingroup\makeatletter\ifx\SetFigFont\undefined%
\gdef\SetFigFont#1#2#3#4#5{%
  \reset@font\fontsize{#1}{#2pt}%
  \fontfamily{#3}\fontseries{#4}\fontshape{#5}%
  \selectfont}%
\fi\endgroup%
{\renewcommand{\dashlinestretch}{30}
\begin{picture}(4224,899)(0,200)
\path(12,864)(612,864)
\path(12,114)(612,114)
\path(87,864)(87,114)
\path(237,864)(237,114)
\drawline(537,114)(537,114)
\path(387,114)(387,118)(387,126)
        (387,139)(387,159)(387,183)
        (387,211)(387,240)(387,269)
        (387,296)(387,322)(387,345)
        (387,365)(387,384)(387,400)
        (387,414)(387,427)(387,439)
        (387,458)(387,475)(387,491)
        (387,507)(387,523)(387,538)
        (387,551)(387,559)(387,563)(387,564)
\path(541,110)(541,111)(540,115)
        (538,125)(534,141)(530,160)
        (525,180)(520,199)(515,216)
        (511,230)(507,243)(503,254)
        (498,265)(493,277)(486,289)
        (478,302)(468,317)(457,333)
        (445,349)(434,364)(425,374)
        (422,379)(421,380)
\path(353,440)(352,441)(350,445)
        (344,453)(338,464)(332,473)
        (328,482)(325,489)(322,496)
        (321,502)(320,508)(319,515)
        (319,522)(320,529)(322,537)
        (325,545)(329,552)(334,560)
        (340,567)(347,575)(355,583)
        (365,591)(377,600)(389,610)
        (401,619)(414,629)(425,638)
        (437,648)(447,656)(457,666)
        (466,675)(475,685)(484,696)
        (492,707)(500,718)(507,729)
        (512,740)(517,750)(521,760)
        (524,770)(527,780)(529,792)
        (531,805)(533,821)(535,837)
        (536,852)(537,862)(537,867)(537,868)
\path(387,673)(387,868)
\path(987,860)(987,110)
\path(912,864)(1512,864)
\path(912,114)(1512,114)
\path(1441,114)(1441,115)(1438,118)
        (1433,126)(1424,139)(1414,154)
        (1402,170)(1391,185)(1381,199)
        (1371,210)(1361,221)(1352,230)
        (1343,239)(1332,247)(1321,256)
        (1308,265)(1293,275)(1276,285)
        (1256,297)(1235,309)(1215,321)
        (1198,330)(1188,336)(1183,339)(1182,339)
\path(1812,114)(2412,114)
\path(1812,864)(2412,864)
\drawline(2337,114)(2337,114)
\path(2337,118)(2337,119)(2335,123)
        (2332,131)(2327,144)(2322,157)
        (2315,170)(2309,182)(2303,193)
        (2295,203)(2287,213)(2279,222)
        (2270,232)(2260,242)(2249,253)
        (2236,264)(2222,275)(2208,287)
        (2194,298)(2179,309)(2165,319)
        (2151,328)(2137,337)(2124,344)
        (2111,351)(2098,358)(2083,366)
        (2066,374)(2047,382)(2027,392)
        (2006,401)(1985,410)(1967,418)
        (1953,424)(1944,427)(1941,429)(1940,429)
\path(1853,463)(1852,464)(1847,469)
        (1839,477)(1832,485)(1828,492)
        (1825,499)(1823,505)(1822,512)
        (1822,519)(1822,527)(1824,535)
        (1827,544)(1831,551)(1836,559)
        (1841,565)(1848,571)(1856,577)
        (1866,583)(1878,590)(1892,597)
        (1907,603)(1924,609)(1941,615)
        (1961,621)(1974,625)(1988,629)
        (2003,633)(2019,637)(2036,641)
        (2054,646)(2072,650)(2091,655)
        (2109,660)(2127,664)(2145,669)
        (2162,674)(2177,678)(2192,683)
        (2205,688)(2218,692)(2236,700)
        (2252,708)(2267,716)(2280,726)
        (2292,735)(2301,745)(2309,754)
        (2316,763)(2321,772)(2325,781)
        (2328,791)(2331,803)(2333,816)
        (2335,832)(2336,849)(2337,863)
        (2337,871)(2337,872)
\path(1887,557)(1887,114)
\path(1887,632)(1887,864)
\path(2041,857)(2041,673)
\path(2037,602)(2037,429)
\path(2033,347)(2033,118)
\path(2187,857)(2187,718)
\path(2191,650)(2191,343)
\path(2187,260)(2187,118)
\path(2712,864)(3312,864)
\path(2712,114)(3312,114)
\path(2787,864)(2787,114)
\path(3237,864)(3237,114)
\path(3091,118)(3091,119)(3090,123)
        (3088,133)(3084,149)(3080,168)
        (3075,188)(3070,207)(3065,224)
        (3061,238)(3057,251)(3053,262)
        (3048,273)(3043,285)(3036,297)
        (3028,310)(3018,325)(3007,341)
        (2995,357)(2984,372)(2975,382)
        (2972,387)(2971,388)
\path(2937,665)(2937,860)
\path(2903,440)(2902,441)(2900,445)
        (2894,453)(2888,464)(2882,473)
        (2878,482)(2875,489)(2872,496)
        (2871,502)(2870,508)(2869,515)
        (2869,522)(2870,529)(2872,537)
        (2875,545)(2879,552)(2884,560)
        (2890,567)(2897,575)(2905,583)
        (2915,591)(2927,600)(2939,610)
        (2951,619)(2964,629)(2975,638)
        (2987,648)(2997,656)(3007,666)
        (3016,675)(3025,685)(3034,696)
        (3042,707)(3050,718)(3057,729)
        (3062,740)(3067,750)(3071,760)
        (3074,770)(3077,780)(3079,792)
        (3081,805)(3083,821)(3085,837)
        (3086,852)(3087,862)(3087,867)(3087,868)
\path(2937,114)(2937,118)(2937,126)
        (2937,139)(2937,159)(2937,183)
        (2937,211)(2937,240)(2937,269)
        (2937,296)(2937,322)(2937,345)
        (2937,365)(2937,384)(2937,400)
        (2937,414)(2937,427)(2937,439)
        (2937,458)(2937,475)(2937,491)
        (2937,507)(2937,523)(2937,538)
        (2937,551)(2937,559)(2937,563)(2937,564)
\path(1103,384)(1102,385)(1097,390)
        (1090,399)(1083,408)(1079,416)
        (1076,425)(1074,433)(1073,442)
        (1072,452)(1072,463)(1072,474)
        (1074,485)(1076,494)(1078,502)
        (1081,510)(1085,518)(1091,525)
        (1098,532)(1106,539)(1115,546)
        (1125,553)(1136,559)(1144,564)
        (1153,569)(1163,575)(1173,581)
        (1185,587)(1197,594)(1209,601)
        (1222,608)(1234,615)(1246,622)
        (1258,629)(1269,636)(1280,643)
        (1291,651)(1302,659)(1314,667)
        (1325,676)(1336,685)(1347,695)
        (1357,704)(1366,713)(1375,723)
        (1382,732)(1388,741)(1395,752)
        (1401,764)(1406,777)(1411,793)
        (1416,811)(1421,831)(1425,848)
        (1428,861)(1429,867)(1429,868)
\path(1137,857)(1137,598)
\path(1137,519)(1137,114)
\path(1287,864)(1287,684)
\path(1287,602)(1287,324)
\path(1287,241)(1287,118)
\drawline(3612,864)(3612,864)
\path(3612,864)(4212,864)
\path(3612,114)(4212,114)
\path(4137,864)(4137,114)
\path(3660,380)(3659,381)(3654,386)
        (3647,395)(3640,404)(3636,412)
        (3633,421)(3631,429)(3630,438)
        (3629,448)(3629,459)(3629,470)
        (3631,481)(3633,490)(3635,498)
        (3638,506)(3642,514)(3648,521)
        (3655,528)(3663,535)(3672,542)
        (3682,549)(3693,555)(3701,560)
        (3710,565)(3720,571)(3730,577)
        (3742,583)(3754,590)(3766,597)
        (3779,604)(3791,611)(3803,618)
        (3815,625)(3826,632)(3837,639)
        (3848,647)(3859,655)(3871,663)
        (3882,672)(3893,681)(3904,691)
        (3914,700)(3923,709)(3932,719)
        (3939,728)(3945,737)(3952,748)
        (3958,760)(3963,773)(3968,789)
        (3973,807)(3978,827)(3982,844)
        (3985,857)(3986,863)(3986,864)
\path(3991,114)(3991,115)(3988,118)
        (3983,126)(3974,139)(3964,154)
        (3952,170)(3941,185)(3931,199)
        (3921,210)(3911,221)(3902,230)
        (3893,239)(3882,247)(3871,256)
        (3858,265)(3843,275)(3826,285)
        (3806,297)(3785,309)(3765,321)
        (3748,330)(3738,336)(3733,339)(3732,339)
\path(3687,523)(3687,118)
\path(3687,857)(3687,598)
\path(3837,856)(3837,676)
\path(3837,591)(3837,313)
\path(3837,241)(3837,118)
\put(50,-30){$0$}
\put(206,-30){$1$}
\put(515,-30){$t$}
\put(365,-30){$s$}
\end{picture}
}

\end{center}
\caption{\label{zopfbild} The braids $\gamma_1,\ldots,\gamma_5$}
\end{figure}
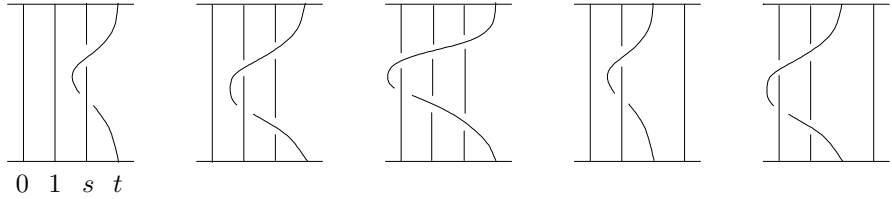

\proof The divisor $D\subset\PP^1_S$ satisfies Assumption
\ref{configass}. Let $\varphi:\pi_1(S,s_0)\to A_4$ be the associated
braiding map. Using standard methods (see e.g.\ \cite{Voelklein01} and
\cite{DettwReiterKatz}), or by staring at Figure \ref{zopfbild}, one
can show that the image of $\varphi$ is generated by the five braids
\[
   \beta_3^2,\;\; \beta_3\beta_2^2\beta_3^{-1},\;\;
   \beta_3\beta_2\beta_1^2\beta_2^{-1}\beta_3^{-1},\;\;
   \beta_2^2,\;\; \beta_2\beta_1^2\beta_2^{-1}.
\]
It is clear that these five braids can be realized as the image under
the map $\varphi$ of generators
$\gamma_1,\ldots,\gamma_5\in\pi_1(S,s_0)$.

 Considering the $\infty$-section as a `tangential base point' for the
fibration $U\to S$ as in \S \ref{mono3}, we obtain a section
$\pi_1(S)\to\pi_1(U)$. We use this section to identify $\pi_1(S)$ with
a subgroup of $\pi_1(U)$. Let $\alpha_1,\ldots,\alpha_5$ be the
standard generators of $\pi_1(U_0)$.  Let $\rho:\pi_1(U)\to K^\times$
denote the representation corresponding to the $G$-cover $f:Y\to X$
and the character $\chi:G\to K^\times$, and $\rho_0:\pi_1(U_0)\to G$
its restriction to the fibre above $s_0$.  Using \eqref{picardeq} one
checks that $\rho_0$ corresponds to the tuple
$\g=(\omega,\omega,\omega,\omega,\omega^2)$, i.e.\ that
$\rho_0(\alpha_i)=g_i$. Also, since the leading coefficient of the
right hand side of \eqref{picardeq} is one, the restriction of $\rho$
to $\pi_1(S)$ is trivial. Hence, by Theorem \ref{etathm}, we have
\[
       \eta(\gamma_i) \;=\; \bar{\Phi}(\g,\varphi(\gamma_i)).
\]
A straightforward computation, using \eqref{locals5eq6} and the
cocycle rule \eqref{locals5eq7}, gives the value of $\eta(\gamma_i)$
(in form of a three-by-three matrix depending on the choice of a basis
of $W_\g$). For this computation, it is convenient to take the classes
of $(1,0,0,0,-\omega^2)$, $(0,1,0,0,-\omega)$ and
$(0,0,1,0,-1)$ as a basis. In order to obtain the $5$ matrices
stated in the theorem, one has to use a different basis, i.e.\ conjugate
with the matrix
\[
     B \;=\; \begin{pmatrix} 
               0          & -\omega-1  & -\omega    \\
               \omega+1   & \omega+1   & \omega+1   \\
               1          & 0          & 0          \\
             \end{pmatrix}.
\]
The claim on the Hermitian form follows from Theorem \ref{poincarethm}
by another straightforward computation.
\Endproof

\begin{rem} 
  Theorem \ref{picardthm} is due to Picard, see \cite{Picard83}, p.\
  125, and \cite{Picard84}, p.\ 181. He obtains exactly the matrices
  given above, but he does not list all of the corresponding braids.
  A similar list as above is obtained in \cite{HolzapfelEuler} using
  different methods.
\end{rem}

\begin{rem}
  It is obvious from Theorem \ref{picardthm} that the Hermitian form
  on $\W$ has signature $(1,2)$ or $(2,1)$, depending on the choice of
  the character $\chi$. This confirms Theorem \ref{signaturethm} in
  this special case. 
\end{rem}




\vspace{4ex}

\begin{minipage}[t]{7cm}
IWR, Universit\"at Heidelberg\\
INF 368\\
69120 Heidelberg\\
michael.dettweiler@iwr.uni-heidelberg.de\\
\end{minipage}
\hfill
\begin{minipage}[t]{4.5cm}
\begin{flushright}
Mathematisches Institut\\ 
Universit\"at Bonn\\
Beringstr. 1, 53115 Bonn\\
wewers@math.uni-bonn.de
\end{flushright}
\end{minipage}

\end{document}